\documentclass[a4paper,12pt]{amsart}

\usepackage{amssymb}
\usepackage{latexsym}
\usepackage{amsthm}

\usepackage[utf8]{inputenc}
\usepackage{color}
\usepackage{multirow}
\usepackage{todonotes}

\newtheorem{thm}{Theorem}
\newtheorem*{thm*}{Theorem}
\newtheorem{cor}[thm]{Corollary}
\newtheorem{lem}[thm]{Lemma}

\newtheorem{prop}[thm]{Proposition}

\theoremstyle{remark}
\newtheorem{rem}[thm]{Remark}
\newtheorem{exmp}[thm]{Example}

\theoremstyle{definition}

\newenvironment{pf}{\par\noindent{\bf Proof.}\enspace\ignorespaces}{\qed\par\par}

\newenvironment{pfm1}{\par\noindent{\bf Proof of Theorem \ref{mainthm}.}\enspace\ignorespaces}{\qed\par\par}

\newenvironment{pfm}{\par\noindent{\bf Proof of Theorem \ref{mainthm} (\ref{mainthm3}).}\enspace\ignorespaces}{\qed\par\par}

\newenvironment{pfm3}{\par\noindent{\bf Proof of Theorem \ref{mainthm4}.}\enspace\ignorespaces}{\qed\par\par}

\newenvironment{pfmu}{\par\noindent{\bf Proof of Theorem \ref{unicitythm}.}\enspace\ignorespaces}{\qed\par\par}

\def\qed{\hfill $\Box$}

\newcommand{\binomial}[2]{\genfrac{(}{)}{0pt}{}{#1}{#2}}

\newcommand{\cha}{\operatorname{char}}

\renewcommand{\div}{\operatorname{div}}

\newcommand{\Pic}{\operatorname{Pic}}

\newcommand{\bQ}{{\mathbb{Q}}}
\newcommand{\bF}{{\mathbb{F}}}
\newcommand{\bR}{{\mathbb{R}}}
\newcommand{\bP}{{\mathbb{P}}}
\newcommand{\bZ}{{\mathbb{Z}}}

\newcommand{\cO}{{\mathcal{O}}}

\newcommand{\cP}{{\mathcal{P}}}
\newcommand{\cL}{{\mathcal{L}}}
\newcommand{\cD}{{\mathcal{D}}}

\DeclareMathOperator{\sep}{sep}

\DeclareMathOperator{\con}{con}

\newcommand{\ggamma}{\overline{\gamma}}
\newcommand{\ogamma}{\gamma^{\sep}}
\newcommand{\cgamma}{\gamma^{\con}}

\DeclareMathOperator{\Gal}{Gal}

\newcommand{\ksep}{{k^{\sep}}}
\newcommand{\kbar}{{\overline{k}}}

\title{Galois descent for the gonality of curves}

\author[Joaquim Ro\'{e} and Xavier Xarles]{Joaquim Ro\'{e} and Xavier Xarles}
\address{Departament de Matem\`atiques\\Universitat Aut\`onoma de
Barcelona\\08193 Bellaterra, Barcelona, Catalonia}
\email{jroe@mat.uab.cat, xarles@mat.uab.cat}

\thanks{The authors were partially supported by the MTM 2013-40680-P (Spanish MICINN grant).
The first author was supported also by the 2014 SGR 114 (Catalan AGAUR
grant).}


\date{\today}

\begin{document}

\begin{abstract}
We determine conditions for the invariance of the gonality under
base extension, depending on the numeric invariants of the curve.
More generally, we study the Galois descent of morphisms of curves to
Brauer-Severi varieties, and
also of rational normal scrolls.
\end{abstract}

\maketitle

\section{Introduction}
Let $C$ be a smooth and projective curve, geometrically connected,
defined over a field $k$ (which, by abuse of language, we call just
a curve). Recall that the gonality $\gamma(C)$ of $C$ over $k$ is an
integer $\gamma\ge 1$ such that there exists a rational non-constant
map $f\colon C\to \bP^1$ of degree $\gamma$ defined over $k$, and
there is no rational map defined over $k$ of degree less than
$\gamma$. Such a rational map $f$ is called a gonal morphism. The
gonality is an important invariant of the curve over $k$, also for
its arithmetic properties (see for example \cite{Fr} and \cite{Po}).
We define the \emph{conic gonality} $\cgamma(C)$ as the minimum
degree of a rational non-constant map $f\colon C\to D$, where $D$ is
a genus 0 curve, all defined over $k$. Finally, we call the gonality
$\gamma(C_{\ksep})$ of $C$ over a separable closure $\ksep$ of $k$
\emph{separable gonality}, and denote it by $\ogamma(C)$, and the
gonality over an algebraic closure $\kbar$ the \emph{geometric
gonality}, and denote it by $\ggamma(C)$. Clearly $$\ggamma(C) \le
\ogamma(C)\le \cgamma(C)\le \gamma(C),$$ since the separably closed
fields are pseudo-algebraically closed, so in particular every conic
over $\ksep$ is isomorphic to $\bP^1$.

In this work we study relations between the gonality $\gamma$, the
conic gonality $\cgamma$ and the separable gonality $\ogamma$ of
$C$, with the aim to obtain sufficient conditions for equalities
between them. These emerge as generalizations of the well known
results for hyperelliptic curves that a curve of genus $g\ge 2$ and
separable gonality $\ogamma=2$ has conic gonality  $\cgamma=2$, and
the result attributed to Mestre \cite{Me}, that a curve of genus
$g\ge 2$ which has even genus and conic gonality $\cgamma=2$, has
gonality $\gamma=2$. Observe that $\gamma=\ogamma$ implies
$\gamma=\cgamma$, but the converse is not true.

Note also that the geometric gonality $\ggamma$ can be smaller than
the separable gonality $\ogamma$: see Example
\ref{genus4sepgonality4} for an example of a genus $4$ curve which
has separable gonality $4$ and geometric gonality $3$.  Observe
however that for genus $\le 3$ they are equal, as well as for
geometric gonality 2.

\medskip
We say that a gonal map $f\colon C\to \bP^1$ is unique if
there is a unique subfield $F$ of the function field $k(C)$ isomorphic to $k(\bP^1)$
and with $[k(C):F]=\gamma$, namely the one determined by the
function $f$.
The main results in this note are the following.

\begin{thm}\label{mainthm} Let $C$ be a curve
with genus $g$ and separable gonality $\ogamma$. Suppose that the
gonal map $f_{\ksep}$ over $\ksep$ is unique. Then
\begin{enumerate}
\item \label{mainunique} $\cgamma=\ogamma$.
\item \label{mainodd}  If the curve $C$ has a $k$-rational divisor of odd degree,
then $\gamma=\cgamma=\ogamma$.
\item\label{maind2} There exists some degree $2$ extension $L/k$ such that
$\gamma(C_L)=\cgamma=\ogamma$.
\item \label{mainthm3} If $\ogamma\equiv
g\pmod 2$, then $\gamma=\ogamma$.
\end{enumerate}
\end{thm}

Brill-Noether theory completely determines the possible gonalities
of a curve of genus $g$ over an algebraically closed field $k$. For
every curve $C$ of genus $g>0$, we have $$2\le \gamma \le
\left\lfloor\frac{g+3}{2}\right\rfloor,$$ and there exists a curve
$C$ of genus $g$ and gonality $\gamma$ for any such number. The
uniqueness hypothesis in Theorem \ref{mainthm} is satisfied in most
cases with non-maximal geometric gonality, at least when the
characteristic of the field is $0$.

\begin{thm*}[Arbarello-Cornalba, {\cite[2.4 and 2.6]{AC}}]
 Let $k$ be algebraically closed and of characteristic $0$, and denote
 $\mathcal{M}^1_{g,d}$ the moduli space of curves of genus $g$
 admitting at least a map of degree $d$ to $\bP^1$, with $g\ge 3$
 and $2\le d< \lfloor(g+3)/2\rfloor$. Then the generic curve in
 $\mathcal{M}^1_{g,d}$ has a unique map of degree $d$ to $\bP^1$.
 More precisely, the locus in $\mathcal{M}^1_{g,d}$
 of curves with two or more such maps
 has codimension at least $g+2-2d\ge 1$.
\end{thm*}

However, it will be much more useful to have effective criteria to decide
whether a given curve has a unique gonal map. We derive such criteria
from Castelnuovo type inequalities like those of \cite{Acc79}.
Recall that a non-constant rational map $f\colon C\to D$ between two
algebraic curves $C$ and $D$ is called simple if there is no smooth
curve $D'$ with maps $f_1 : C \to D'$, and $f_2 : D'\to  D$ such
that $\deg(f_i) \ge 2$, for $i = 1, 2$, and $f = f_2 \circ f_1$.
Equivalently, if the corresponding extension of function fields
$k(C)/k(D)$ is simple as extension. For example, if the degree of
$f$ is a prime number, then $f$ is simple.
Under a simplicity hypothesis and
a suitable bound on the separable gonality, we can prove that
uniqueness holds:

\begin{thm}\label{unicitythm}
 Let $C$ be a curve with genus $g$ and separable gonality $\ogamma$.
  The gonal map over $\ksep$ is unique if it is simple
 (in particular this is always the case if  $\ogamma$ is prime)
  and $(\ogamma-1)^2<g$.
\end{thm}

From now on we say that a curve of genus $g$ has low gonality if
$(\ogamma-1)^2<g$.

\medskip

 Note also that the class of \emph{goneric} curves introduced in
 \cite{SSW}, have unique gonal maps over $\kbar$.
Gonericity is a condition expressed in terms of the gonality and the Betti numbers
of a minimal resolution of the ideal of the curve in its canonical embedding,
 and there exist efficient algorithms to decide whether
 a given curve is goneric. In fact, proposition 3 in \cite{SSW}
 shows that if $C$ is goneric, then $W_{\ggamma}^1$ is a single
 reduced point, which shows that
 in this case (since $\ksep$ is pseudo-algebraically closed)
 $\ogamma=\ggamma$ as well.

\medskip

We give two proofs for theorem \ref{mainthm}. The first one, given
in section \ref{sec:brauer-severi}, is based on the theory of
Brauer-Severi varieties, and it naturally leads to analogous results
for maps from $C$ to $\bP^r$:  in this general case, uniqueness over
$\ksep$ implies descent to a map to a Brauer-Severi variety (see
Theorem \ref{descunique} for a precise statement).

A curve in projective space over an algebraically closed field is
called reflexive if the composition of its Gauss map with that of
its dual is an isomorphism, which is always the case in
characteristic zero, and in characteristic $p>2$ is equivalent to
the intersection multiplicity with its tangent line at a general
point being 2 \cite[3.5]{HK85}; we say that $C$ is reflexive if
$C_{\kbar}$ is reflexive. Then, we can give a generalized version of
Theorem \ref{mainthm} (\ref{mainthm3}) to maps to $\bP^r$ for $r>1$;
in the case $r=2$ we obtain the following result.

\begin{thm}\label{mainthm4} Given a curve $C$ be a curve defined over $k$,
denote $\gamma_2$ (resp. $\ogamma_2$) the smallest degree of a plane
model of $C$ (resp. of $C_{\ksep}$).
\begin{enumerate}
\item  Suppose that the
corresponding $g^2_{\ogamma_2}$ on $C_{\ksep}$ is unique. If
$\ogamma_2\not\equiv 0\pmod 3$, then $\gamma_2=\ogamma_2$.
\item If the plane model of $C_{\ksep}$ is reflexive
(for instance, if $\cha k=0$) and $g> \left\lfloor
\frac{(\ogamma_2)^2-3\ogamma_2+3}3 \right\rfloor$ then the
corresponding $g^2_{\ogamma_2}$ is unique.
\end{enumerate}
\end{thm}

In section \ref{sec:rns} a second proof of Theorem \ref{mainthm}
(\ref{mainthm3}) is obtained, which gives additional information, at
least for geometrically trigonal curves (see theorem
\ref{theoremtrig}). This proof is based on the study of Galois
descent for rational normal scrolls, which may be of independent
interest.

\smallskip
\noindent\textbf{Acknowledgements.} We thank Ciro Ciliberto for a
careful reading of the first version of this paper and for pointing
us to relevant results on Castelnuovo inequalities, and Miguel \'{A}ngel
Barja, Andr\'{e} Hirschowitz and Bjorn Poonen for helpful conversations
and suggestions. We thank the referee of the paper for some
corrections, as well as for suggesting to study the geometric gonality
and its relation to the separable gonality.

\section{Galois descent and uniqueness}
\label{sec:descent} Let $C$ be a curve defined over a field $k$, and
fix a separable and an algebraic closure $k\subset\ksep\subset \kbar$
for the whole paper.

Given any divisor $D$ defined
over $k$, denote by $\cL(D):=H^0(C,\cO(D))$ the $k$-vector space of
meromorphic functions $f$ on $C$ such that $\div(f)+D \ge 0$ is
effective. Denote also by $|D|$ the set of effective divisors
linearly equivalent to $D$. Then there is a canonical bijection
between $\bP(\cL(D))$ and $|D|$ determined by mapping $f\in \cL(D)$
to the divisor $\div(f)+D$.

Recall that an $r$-dimensional linear series $\cD$ over $\ksep$ is
the family of divisors given by a vector subspace $V$ in $\cL(D)$,
for some divisor $D$, and $\cD$ is a $g^r_d$ if $\deg D = d$ and
$\dim V = r+1$. The linear series is complete if $V=\cL(D)$. Recall
also that a base-point-free linear series $\cD$ determines a
$\ksep$-morphism $\phi_{\cD}:C\to \bP(V^*)\cong \bP^r$, which maps a
point $P\in C(\ksep)$ to the point corresponding to the hyperplane
$$\phi_{\cD}(P):=\{s\in V \ | \ s(P)=0\}.$$
A base-point-free linear series $g^r_n$ is called
simple if the map $\phi$ it determines
is simple in a suitable sense, namely it can not be factored
as $\phi=\phi'\circ f$ with $f : C \to C'$ and $\phi' : C'\to  \bP^r$
for some curve $C'$, with $\deg (f),\deg(\phi') \ge 2$
(where $\deg(\phi')=n'$ is the degree of the linear series $g^r_{n'}$
induced on $C'$, i.e., $n=\deg(f_1)\cdot n'$). For
$r=1$ this is equivalent to the map of curves
$C\to \bP^1$ being simple in the sense above, whereas
for $r>1$ it is equivalent to $\phi$ being birational
onto its image.

We say that a complete linear series $\cD$ is Galois invariant if
for any $\sigma\in \Gal(\ksep/k)$, the divisor $D^{\sigma}$ is
linearly equivalent to $D$. In this case, any $E\in |D|$ is also in
$|D^{\sigma}|$, hence one gets a natural action of $\Gal(\ksep/k)$
in $\bP(\cL(D))$. We say that a linear series $\cD$ given by
$V\subset \cL(D)$ is Galois invariant if the corresponding complete
linear series and the subspace $\bP(V)\subset\bP(\cL(D))$ are Galois
invariant.

It is surely well known to the experts that morphisms
$\pi:C\rightarrow C'$ of degree $d$, with $C'$ a genus zero curve,
correspond to $g^1_d$'s on $C_{\ksep}$ invariant under the Galois
action. We want to extend this result to higher dimensions of the
target space.

Recall that a Brauer-Severi variety of dimension $r$ over $k$ is
a smooth projective variety $\cP$ such that $\cP\otimes_k \ksep\cong
\bP^r_{\ksep}$. Hence the Brauer-Severi varieties $\cP$ of dimension
$1$ are the curves of genus $0$.

\begin{lem}\label{pencilconic}
Let $X$ be a variety defined over $k$. Then the set of morphisms to
some $r$-dimensional Brauer-Severi variety, modulo automorphisms,
corresponds bijectively to the base-point free $r$-dimensional
linear series over $\ksep$ invariant under the Galois action of the
absolute Galois group $\Gal(\ksep/k)$.
\end{lem}

\begin{pf} That morphisms to Brauer-Severi varieties give
$r$-dimensional linear series over $\ksep$ invariant under the
Galois action is clear.

To show bijectivity, first of all observe that, if a divisor $D$ is
invariant by the action of  $\Gal(\ksep/k)$, and the associated
linear series $\cD$ is base-point free, then it determines a
$k$-defined map $\phi_{\cD}:X_k\to \bP(V^*)_k\cong \bP_k^r$.

Now, if the linear series $\cD$ is Galois invariant, one gets an
action of $\Gal(\ksep/k)$ in $\bP(V)$  by automorphisms; hence also
a dual action on $\bP(V^*)$. Both actions determine Brauer-Severi
varieties $\cP$ and $\cP^*$ by the classical theory of Galois
descent, which become split on the field of definition of the
corresponding divisor $D$.

Finally, we only need to show that the map $\phi_{\cD}$ commutes
with the action of $\Gal(\ksep/k)$. But
$$\phi_{\cD}(P^{\sigma}):=\{s\in V \ | \
s(P^{\sigma})=0\}=\{s^{\sigma^{-1}}\in V \ | \
s^{\sigma^{-1}}(P)=0\}=\phi_{\cD}(P)^{\sigma},$$ since the dual
action of $\sigma$ in $\bP(V^*)$ sends a point corresponding to a
subspace $W\subset V$ to the subspace $W^{\sigma^{-1}}$.
\end{pf}

As a corollary, one immediately obtains:
\begin{thm}\label{descunique}
Let $C$ be a (smooth projective) curve defined over $k$. Suppose
that for a fixed $r$ and $d$ there is only one $g^r_d$, giving a
morphism $f:C_{\ksep}\to \bP^r_{\ksep}$. Then there exists a
Brauer-Severi variety $\cP$ defined over $k$ together with a
$k$-morphism $g:C\to \cP$ such that $g\otimes_{k}\ksep:C_{\ksep}\to
\cP_{\ksep}\cong\bP^r_{\ksep}$ is equal to $f$.
\end{thm}

Determining the uniqueness of $g^r_d$'s for $r$ and $d$ small
relative to the genus is a classic problem (see  \cite{Acc79},
\cite{Cil83}, \cite{CL}, \cite{CK90}, and references therein). In
particular it is widely known that the $g^2_d$ of a smooth plane
curve of degree $d>3$ is unique, but the same is true for ``mild''
singularities, i.e., if $d$ is small enough compared to the genus.

We approach uniqueness of $g^r_d$'s by the classical ``Castelnuovo
method'', see \cite{Cas}, \cite{Acc79}, \cite{Cil83},
\cite[IV.6]{Har}. The idea is to estimate the dimension of sums of
linear series, by counting conditions. To begin with, if $\mathcal
D_1$ and $\mathcal D_2$ are linear series, the sum $\mathcal
D_1+\mathcal D_2$ is the minimal linear series containing all
divisors $D_1+D_2$ with $D_i\in \mathcal D_i$. If $\mathcal D_i$ is
given by the linear subspace $V_i\subset H^0(C,\cO(D_i))$, then
$\mathcal D_1+\mathcal D_2$ is given by the image of the map
$$\phi:V_1\otimes V_2\rightarrow H^0(C,\cO(D_1+D_2))$$
determined by $\phi(s_1\otimes s_2)=s_1s_2$. On the other hand, if
$\mathcal D$ is a linear series given by the
linear subspace $V\subset H^0(C,\cO(D))$, and $D'$ is an arbitrary
effective divisor, one puts
$$ V(-D')=\{s\in V \,|\,\div(s)\ge D'\}.$$
The linear series determined by $V(-D')$ has $D'$
as a fixed divisor; subtracting $D'$ from it one obtains
a series denoted $\mathcal D -D'$. The \emph{number of conditions}
imposed by $D'$ on $\mathcal D$ is $\dim \mathcal D -\dim (\mathcal D-D')$,
or equivalently $\dim V- \dim (V(-D'))$.

Recall that a curve $C$ in projective space $\bP^r_k$ is called
reflexive if the composition of the Gauss map of $C_{\kbar}$ with
that of its dual is an isomorphism, and it is called strange if
there is a point of $\bP^r_{\kbar}$ that belongs to every tangent
line of $C_{\kbar}$. Observe that reflexive curves are not strange.
On the other hand, the only nonsingular strange curves over an
algebraically closed field are lines, and conics in characteristic 2
\cite[IV.3.9]{Har}, but the plane smooth curve
$x^{p+1}+y^{p+1}+z^{p+1}=0$ is not reflexive in characteristic $p$
(see \cite{HK85}, \cite[2.3]{Rat} for this and other examples) so
there are indeed curves that are neither reflexive nor strange.

\begin{lem} \label{unifpos}
Let $k$ be an algebraically closed field, and $C$ a curve defined over $k$.
 Suppose $\mathcal D_1$ and $\mathcal D_2$ are two different
 base-point-free simple $g^r_n$'s on $C$, $n>r\ge 2$,
 and let $f:C\rightarrow \bP^r_{k}$
 be the map induced by $\mathcal D_1$. Assume that at least one
 of the following is true:
\begin{enumerate}
\item $f(C)$ is reflexive,
\item $f(C)$ is not strange
(for instance, $f(C)$ is nonsingular) and $r\ge 4$.
\end{enumerate}
Then a general divisor $D \in \mathcal D_1$ is made up of
$n$ distinct points, and every subset of $r$ points in $D$
imposes $r$ conditions to  $\mathcal D_1$, whereas
every subset of $r+1$ points in $D$
imposes $r+1$ conditions to  $\mathcal D_2$.
\end{lem}

Recall that a statement claimed for \emph{a general divisor in $\cD$}
is meant to hold for all
divisors in a nonempty Zariski-open subset of $\cD \cong \bP^r$.

\begin{proof}
 That over an algebraically
  closed field a general divisor is made up of
 distinct points  is well known (\cite[IV, Exercise 3.9]{Har}, \cite[Lemma 3.11.2]{St}).

Then, a general divisor $D\in \cD_1$ imposes $r$ conditions to
 $\cD_1$ (i.e., $D$ is the unique divisor in $\cD_1$ containing $D$)
 and $r+1$ conditions to
 $\cD_2$ (otherwise $D\in \cD_1\cap \cD_2$, but if this happens for
 general $D\in \cD_1$, the two linear series must coincide).
 Under the hypotheses, the monodromy group of $f$
 contains $A_n$ \cite[2.2, 2.5]{Rat}. Then the same proof as in \cite[1.8]{Rat}
 gives the result.
 \end{proof}

Following Accola, we set $$R(l;r)=l(l+1)r/2-l(l-1)/2$$  and
$$R(l_1,l_2;r)=R(l_1,r)+R(l_2,r)+l_1l_2r.$$ Castelnuovo's method (see
\cite{Acc79}, \cite[Lemma 1.3]{Cil83}) yields the following:

\begin{lem}[Accola, {\cite[4.2]{Acc79}}]\label{accolacharp}
  Suppose $\mathcal D_1$, and $\mathcal D_2$ are two different
 simple $g^r_n$'s without
 fixed points on $C$,  and assume that
 a general divisor $D_i \in \mathcal D_i$ is made up of
$n$ distinct points, and
 \begin{enumerate}
 \item every subset of $r$ points in $D_i$
imposes $r$ conditions to  $\mathcal D_i$,
\item every subset of $r+1$ points in $D_i$
imposes $r+1$ conditions to  $\mathcal D_j$, $j\ne i$.
 \end{enumerate}
Then $\dim (l_1\cD_1+l_2\cD_2)\ge R(l_1,l_2,r)$ for all non-negative
integers $l_1,l_2$ satisfying $(l_1+l_2)r+l_1-1\le d$.
 \end{lem}

\begin{thm}\label{accolaunique}
 Let $C$ be a curve of genus $g$ over an
 algebraically closed field $k$, with a simple linear
 series $g^r_d$, giving a morphism
 $f:C\to \bP^r_{k}$. Assume that at least one of the
 following is true:
\begin{enumerate}
\item $k$ has characteristic zero,
\item $r=1$,
\item $f(C)$ is reflexive, and $r\ge 2$,
\item $f(C)$ is not strange
(for instance, $f(C)$ is nonsingular) and $r\ge 4$.
\end{enumerate}
 Write $d = m(2r - 1) + q$ where $q$ is the residue of
 $d$ modulo $(2r- 1)$ so that $-r + 2 \le q \le r$. Let $v=1$ if $q\le 1$,
 $v=0$ otherwise. If
\begin{equation}
\label{eq:accola}
  g >m^2(2r-1)+m(2q-1-r)-v(q-1)
\end{equation}
 then the given series is the unique simple $g^r_d$ on $C$.
\end{thm}
\begin{rem}
 For $r=1$, $f(C_{\ksep})=\bP^1_{\ksep}$, $q=v=1$,
 $m=d-1$ and the inequality \eqref{eq:accola} reads simply $g>(d-1)^2$. In this
 case, the result of Proposition \ref{accolaunique} was already known
 to Riemann \cite{Rie57} (for $k=\mathbb C$).
 For $r=2$, the inequality \eqref{eq:accola} is equivalent to
 $$
 g> \left\lfloor \frac{d^2-3d+3}3 \right\rfloor.
 $$
\end{rem}

\begin{proof}
For $r=1$, the result follows from the so-called
``Castelnuovo-Severi'' inequality, see for instance
\cite[III.11.3]{St}. So assume $r\ge 2$.

 When $k$ has characteristic zero, R. Accola has shown
 that the existence of two distinct simple $g^r_d$'s
 contradicts the inequality \eqref{eq:accola}, in
 \cite[Theorem 4.3]{Acc79}. The proof relies on a uniform
 position lemma \cite[4.1]{Acc79}, which needs characteristic
 zero, to show that the number of conditions imposed by divisors
 satisfies the hypotheses of lemma \ref{accolacharp}.
 In our case they are satisfied thanks to  \cite[1.8]{Rat}
 and lemma \ref{unifpos}.
 The rest of Accola's argument consists in matching the dimension
 estimate of lemma \ref{accolacharp} with Clifford's inequality for
 special divisors. This does not depend on the characteristic,
 so the result follows.
\end{proof}

We will use the preceding results, which are stated for
algebraically closed fields, in the proof of Theorem
\ref{unicitythm}; therefore, we need to consider the inseparable
base change $\kbar / \ksep$. It is probably
well known to the experts that under the key assumption made
throughout the paper that $C$ is a smooth curve, the relevant
phenomena are all stable under this base change:

\begin{lem}\label{inseparable}
Let $C$ be a (smooth, projetive, geometrically connected) curve of
genus $g$ over an arbitrary field $k$, and fix an algebraic closure
$\kbar$ of $k$. Let $f:C\to D$ a morphism where $D$ is a (non
necessarily smooth) projective curve defined over $k$, and for every
algebraic extension $L/k$, denote $ f_L:C_{L}\rightarrow D_{L}$, its
base change.
\begin{enumerate}
\item For every algebraic extension $L/k$, the field of rational
functions of $C_L$ is separable over $L$, of genus $g$.
\item For every algebraic extension $L/k$, $\deg f_{L}=\deg f$.
\item If  $f_{\ksep}$ is simple, then $f_{\kbar}$ is simple.
\item If a $g^r_d$ is simple on $C_{\ksep}$, then it is simple
on $C_{\kbar}$.
\end{enumerate}
\end{lem}

\begin{proof}
If $k$ is perfect then all claims are well known (see
\cite[chapter III]{St}) so assume $k$ is an imperfect field of characteristic $p$.

  Since $C$ is smooth, $C_L$ is smooth over $L$.
In particular the field $k(C_L)$ is formally smooth over $L$,
and therefore separable.
Moreover, genus does not change for function fields of
smooth curves, by \cite[\S 3]{Ros52}.

For the second claim, we first prove that
the fields $L$ and $k(C)$ are linearly disjoint over $k$.
By the transitivity of linear disjointness \cite[VIII,3.1]{Lan02}
it is enough to consider the cases
that $L/k$ is separable or purely inseparable.
In the separable case, as both $k(C)$ and $L$ are separable,
the proof of \cite[III.6.1]{St} shows that $k(C)$ and $L$ are
linearly disjoint over $k$. In the inseparable case,
because $k(C)$ is separable over $k$,
it is linearly disjoint with $k^{p^{-\infty}}$
(MacLane's criterion, \cite[VIII,4.1]{Lan02})
and hence with $L\subset k^{p^{-\infty}}$ over $k$.
Now $L$ and $k(C)$ being linearly disjoint over $k$
implies that $k(D_L)=Lk(D)$ and $k(C)$ are
linearly disjoint over $k(D)$. This gives the second claim.

For the third claim, assume by way of contradiction that $f$ is
simple and there is a nontrivial intermediate field $F$,
$k(C_{\kbar})\supsetneq F \supsetneq k(D_{\kbar})$. Then $F$ is
purely inseparable over $F \cap k(C_{\ksep})$, and by the simplicity
of $f$, $F \cap k(C_{\ksep})=k(D_{\ksep})$. So $F$ is purely
inseparable over $k(D_{\ksep})$ and therefore over $k(D_{\kbar})$ as
well. This implies that there is an element $t\in k(D_{\kbar})$
which is not a $p$th power, with $t^{1/p}\in F\subsetneq k(C_{\kbar})$.
On the other hand, let $x\in k(C_{\ksep})\setminus
k(D_{\ksep})$ and let $P(X)$ be its minimal polynomial over
$k(D_{\ksep})$. By the simplicity of $f$,
$k(C_{\ksep})=k(D_{\ksep})(x)$ and $\deg P=\deg f=n$. Moreover
$k(C_{\kbar})=k(D_{\kbar})(x)$, and by the second claim, $P$ is also
the minimal polynomial of $x$ over $k(D_{\kbar})$. In particular,
since $k(D_{\kbar})(t^{1/p})\subsetneq k(C_{\kbar})$, $p$ is a
proper divisor of $n$. Now let $Q(X)$ be the minimal polynomial of
$x$ over $k(D_{\kbar})(t^{1/p})$. Considering the degrees of the
extensions, one clearly has $\deg Q=n/p$. But then $Q^p\in k(D_{\kbar})
[X^p]$ is a monic polynomial of degree $n$ in $X$ vanishing at
$x$, i.e., $Q^p=P$. This means that $P$ only involves $p$th powers
of $X$, i.e., $x^p$ is an element in $k(C_{\ksep})$ whose minimal
polynomial has degree $n/p<p$, contradicting the simplicity of $f$.

Finally, a linear series $g^r_d$ with $r\ge 2$ is simple if and only if
the map $f:C \rightarrow \bP^r$ it defines is birational onto
its image. Therefore, by the second claim simplicity does not
change under base field extension. The case $r=1$ has been dealt
with in the third claim.
\end{proof}

\begin{rem}
  Counterexamples to lemma \ref{inseparable}
when $C$ is not smooth do exist. For instance one can consult
articles on ``genus change under inseparable extensions'', starting
from classical Tate's and Rosenlicht's papers \cite{Tat52},
\cite{Ros52}.
\end{rem}

\begin{pfmu}
%
Let $f$ be a gonal map over $\ksep$, which by assumption is simple.
By lemma \ref{inseparable} its base change $f_{\kbar}:C_{\kbar}
\rightarrow \bP^1_{\kbar}$ is simple as well. By theorem
\ref{accolaunique} in the case $r=1$, the unique simple
$g^1_{\ogamma}$ on $C_{\kbar}$ is then the one determined by
$f_{\kbar}$. A fortiori, the unique simple  $g^1_{\ogamma}$ on
$C_{\ksep}$ is the one determined by $f$.
\end{pfmu}

\begin{exmp} It is well known that there exist curves  with genus
$g=(\ogamma-1)^2$ and with more than one simple gonal map. For
example, a curve $C$ embedded as a smooth curve of type $(\ogamma,
\ogamma)$ in the smooth quadric surface $Q=\bP^1 {\times} \bP^1 \subset
\bP^3$. In this case $C$ has gonality $\ogamma$ and it has exactly
two $g^1_{\ogamma}$: the ones induced by the projections of $Q$ onto
one of its factors (see for example \cite{Mar}, Theorem 1).
\end{exmp}

\begin{rem}
 If $r=2$ and
 $g> \left\lfloor \frac{d^2-3d+3}3 \right\rfloor$,
 so that theorem \ref{accolaunique} holds, every
 $g^1_{e}$ on $C$ with $e<d$ is cut out on $f(C_{\ksep})\subset \bP^2_{\ksep}$
 by a pencil of lines, and there is no $g^2_{e}$ with $e<d$.
 So the separable gonality is $d-m$ where $m$ is the maximal multiplicity
 of a singular point of $\Gamma$.
\end{rem}

We include the reference to a last uniqueness result, due to
Cilliberto and Lazarsfeld \cite{CL}, which applies to $r=3$ in the
case of complete intersections.

\begin{thm*}[Ciliberto--Lazarsfeld] Let $C\subset \bP^3$ be a smooth
curve defined over a field of characteristic $0$, which is the
complete intersection of two surfaces of degree $h$ and $h'$, both
bigger than $4$. Then any simple $g^s_{m}$, with $m\le hh'$ and
$s\ge 2$ is unique (and, in particular, the canonical $g^3_{hh'}$).
\end{thm*}

\section{Splitting Brauer-Severi Varities}
\label{sec:brauer-severi}

Recall that to any Brauer-Severi variety $\cP$ of dimension $n$ over
a field $k$ one can assign canonically a central simple algebra of
rank $(n+1)^2$ over $k$. Its class $[x]$ in the Brauer group
verifies that in the exact sequence
$$ \Pic(\cP) \to \Pic(\cP\otimes_k{\ksep})\cong \bZ \to Br(k) $$
the last map sends $1$ to $[x]$, and hence the image of some
generator of $\Pic(\cP)$ is equal to $m$, where $m$ is the order of
$[x]$. Hence $m$ divides $n+1$ since $[x]$ has order dividing $n+1$.
We say that $\cP$ is split over $k$ if $\cP\cong \bP^n$ already over
$k$. Then $\cP$ is split if and only if $[x]=0$, i.e. its order is
equal to $1$.

The following result summarizes some properties of $\cP$, some of
which are well known.

\begin{thm}
\label{BSsplitting}
Let $\cP$ be a Brauer-Severi variety of dimension $n$
over a field $k$. Then
\begin{enumerate}
\item \label{BS1} If $\cP$ contains a hypersurface of degree
coprime with $n+1$, then $\cP$ is split.
\item \label{BS2} If there is a Galois invariant element in $\bZ[\cP(\ksep)]$
(e.g. it $\cP$ has a $k$-rational point) whose degree is coprime
with $n+1$, then $\cP$ is split.
\item \label{BS3} There exists an immersion $\cP \to \bP^N$, where $N=\binomial{2n+1}{n}-1 $, as a smooth subvariety of degree $(n+1)^n$.
\item \label{BS4} There exists a finite map $f:\cP \to \bP^n$ of degree $(n+1)^n$.
\item \label{BS5} There exists an extension $L/k$ of degree dividing $n+1$ such that
$\cP\otimes_k L\cong \bP^n_L$.
\end{enumerate}
\end{thm}

\begin{pf} The assertion \eqref{BS1} is clear since the hypersurface determines
an element in $\Pic(\cP)$ whose image in $\Pic(\cP\otimes_k{\ksep})$
is equal to the degree of the hypersurface.

The next result \eqref{BS2} is a generalization of a result by
Ch\^{a}telet in \cite{Cha} (see also \cite{Ar}), who showed the case
of rational points. We will use the dual Brauer-Severi variety
$\widehat{\cP}$. It is a Brauer-Severi variety together with an
inclusion reversing correspondence between twisted linear
subvarieties of dimension $d - 1$ of $\cP$ and those of codimension
$d$ in $\widehat{\cP}$. Now, a Galois invariant element in
$\bZ[\cP(\ksep)]$ determines an element in $\Pic(\widehat{\cP})$,
whose degree is equal to the degree of the formal sum, since the
degree of a linear hypersurface is 1. The result is deduced then
from \eqref{BS1}.

Result \eqref{BS3} is well known: in fact, the immersion is given by
the anticanonical sheaf, which is always defined over $k$. It is
known that the anticanonical sheaf in $\bP^n$ is equal to
$\cO(n+1)$, and it gives the $(n+1)$-tuple Veronese embedding in
$\bP^N$.

Now, choosing $n-1$ sufficiently general hyperplanes in $\bP^N$, we
can find some whose intersection, which is a linear subvariety of
dimension $N-n-1$, does not intersect the image of $\cP$ in $\bP^N$.
Projecting to a complementary linear subvariety of dimension $n$ we
get the desired finite morphism.

The last result is due to F. Ch\^{a}telet in his thesis, and it is a
consequence of the main classical results on central simple
algebras. If $A$ denotes a central simple algebra associated to
$\cP$, then $\cP$ splits over an extension $L/K$ if and only if $A$
does. But $A$ always splits over a maximal commutative subfield,
which has degree over $K$ equal to the index of $A$, which is the
square root of the dimension of the associated division algebra,
which clearly divides $n+1$.
\end{pf}

\begin{cor}\label{Conicsplitting} Let $C$ be a curve defined over a field $k$,
with genus $g$, gonality $\gamma$ and conic gonality $\cgamma$. Then
\begin{enumerate}
\item \label{conicbound} $\cgamma\le\gamma\le 2\cgamma$.
\item \label{conicodd} If the curve $C$ has a $k$-rational divisor of odd degree,
then $\gamma=\cgamma$.
\item \label{coniceven} If $\cgamma\ne\gamma$, then $\gamma$ is even.
\item \label{conicext} There exists some degree $2$ Galois extension $L/k$ such that
$\gamma(C_L)=\cgamma$.
\end{enumerate}
\end{cor}
\begin{pf}
 \eqref{conicbound} (resp. \eqref{conicodd}) follows from \eqref{BS4}
 (resp. \eqref{BS2}) in theorem \ref{BSsplitting},  in the case
 $n=1$. The point \eqref{conicext} follows from \eqref{BS5}, except the fact
 that we can take the degree 2 extension to be Galois. This follows from the
 fact that any conic has a point in a separable extension of degree
 2. If the characteristic of the field is not $2$, this is clear. If
 it is 2, and it has no $k$-rational point, then the conic can be described in $\bP^2$ by an equation of the form
$ax^2+by^2+cz^2+xz+yz=0$ for some $a$, $b$ and $c\in k^*$. Then the
points which intersect the line $x=0$ determine the desired
extension. Finally, \eqref{coniceven} is immediate from
\eqref{conicodd}.
\end{pf}

Observe that it is not true that the gonality is always the conic
gonality or its double, as the following example shows.

\begin{exmp} The genus 4 curve over $\bQ$ (or even over $\bR$) given in canonical form as
the intersection in $\bP^3$ of the quadric $x^2+y^2+z^2=0$ with the
cubic $x^3+y^3+t^3=0$, has conic gonality $3$ (with unique gonal map
given by the projection map to the conic $x^2+y^2+z^2=0$), and
gonality $4$ (with gonal map given by the natural projection to the
cubic $x^3+y^3+z^3=0$ followed by the degree two map determined by a
rational point of the cubic (e.g. $[1:-1:0]$)).
\end{exmp}

\begin{exmp}\label{genus4sepgonality4} Let $k=\bF_2(s)$ be the
field of rational functions over the finite field $\bF_2$. The
genus 4 curve given in canonical form as the intersection
in $\bP^3$ of the quadric $xy+z^2+st^2=0$ with the cubic
$x^3+y^3+t^3=0$, has geometric gonality $3$, and separable gonality
$4$ (with a gonal map given as in the previous example). In fact the
gonal maps over $\kbar$ are given by the two rulings of the quadric,
which are not defined over $\ksep$.
\end{exmp}

\begin{pfm1}
Part (1) is corollary \ref{descunique} in the case $r=1$.
Parts (2) and (3) are then consequence of the first part and
corollary \ref{Conicsplitting}.
\end{pfm1}

\begin{prop}\label{mainthm2} Let $C$ be a curve defined over a field $k$
with genus $g$ and conic gonality $\cgamma$. Assume that over a
separable closure $\ksep$, the $g^1_{\cgamma}$ associated to a
conic-gonal map is complete. If $\cgamma\equiv g\pmod 2$, then
$\gamma=\cgamma$.
\end{prop}

\begin{pf}
 Let $f:C\rightarrow D$ be a map of degree $\cgamma$
to a curve of genus zero as in the claim. If $D\cong \bP^1$ there is
nothing to prove, so assume $D$ is a conic. By corollary
\ref{Conicsplitting} (4), there is a degree 2 Galois extension $L$
of $k$ over which $D$ becomes split. Consider the Galois-invariant
$g^1_{\cgamma}$ on $C_{L}$, which by hypothesis is complete; call
$E$ a divisor in this $g^1_{\cgamma}$. Note that, although $E$ is
only defined over $L$, the linear series $|E|$ is Galois-invariant.

Now if $K$ is a canonical divisor
the complete linear series $|K-E|$
is Galois-invariant because $K$ and $|E|$ are Galois-invariant.
By Riemann-Roch, and because the $g^1_{\cgamma}$
is complete, $\dim |K-E|=g-\cgamma=n$,
and we get a morphism $C\rightarrow \mathcal P$
to a Brauer-Severi variety of (even) dimension $n$.

Since $K-E$ is a divisor defined over $L$,
$\mathcal P\otimes_k L \cong \bP^n_L$.
Now if $H$ is a hyperplane defined over $L$, and
$H^\sigma$ is its conjugate by $\Gal(L/k)$,
$H+H^\sigma$ is a degree 2 hypersurface
defined over $k$. So ${\mathcal P}$ contains a hypersurface of degree 2,
which is coprime with $n+1$, and by theorem \ref{BSsplitting},
${\mathcal P}\cong \bP^n$. This means that $K-E$
(and hence $E$) is linearly equivalent to a divisor
defined over $k$, so $D\cong \bP^1$
and we are done.
\end{pf}

\begin{rem}
 We have seen in the course of the proof that, if
over a separable closure $\ksep$, the $g^1_{\cgamma}$ associated to
a conic-gonal map is complete, then $\cgamma = g-\dim |K-E|\le g+1$.
\end{rem}

\begin{pfm}
  By theorem \ref{mainthm}, $\cgamma=\ogamma$, and
  because gonal series over a separably closed field are
  always complete, proposition \ref{mainthm2} applies, so $\gamma=\cgamma$.
\end{pfm}
\begin{pfm3}
  By corollary \ref{descunique}, if the $g^2_{\ogamma_2}$ is unique
  there is a morphism $g:C\to \cP$, with $\cP$ a Brauer-Severi
  variety of dimension $n=2$, whose base change to $\ksep$ is the
  generically injective morphism of lowest degree.
  $g(C)$ is a divisor of degree $\ogamma_2$, coprime with $n+1=3$,
  so by theorem \ref{BSsplitting} \eqref{BS1}, $\cP\cong \bP^2$.

 The second claim is a direct application of
theorem \ref{accolaunique} in the case $r=2$,
taking into account lemma \ref{inseparable}.
\end{pfm3}

\begin{rem}\label{planedescent}
 Let $C$ be a curve such that $C_{\ksep}$ has a plane model
 $\Gamma$ of degree $d$ not divisible by $3$,
 such that all
singularities of $\Gamma$ are nodes or ordinary cusps,
 and assume that either
 $(d-2)(d-3)<2g-2$ or
 $d\ge 8$ and $(d-3)(d-4)<2g-10$ .
Then $\gamma_2=\ggamma_2=d$ and $\Gamma$ is defined over $k$. In
particular, the smooth plane curves over $\ksep$ of degree $d>3$ and not multiple
of $3$ are already plane curves over $k$.

Indeed, by Coppens-Kato \cite[Theorem 2.4]{CK90} and \cite{CK91},
the $g^2_d$ on $C_{\ksep}$
 is unique, and $\ogamma_2=d$. Then by theorem \ref{mainthm4},
 $\gamma_2=\ggamma_2$ and $\Gamma$ is defined over $k$.
 Note that the nature of the proofs in the Coppens-Kato papers
 is independent of the characteristic, and hold without the
 reflexivity hypothesis. In contrast, the inequalities on the
genus are more restrictive than in theorem \ref{unicitythm}.
\end{rem}

\begin{exmp} A smooth plane curve of degree $3$ has genus 1, and
over $\ksep$ any such curve has such a plane model and gonality 2.
It is well known that there are genus one curves over $\bQ$ (or over
any number field) with gonality $d$ for any $d\ge 2$ (see
\cite{Clark}, and use that the index and the gonality are equal if
the index is larger than one by Riemann-Roch). Hence, if the
gonality is larger than $3$, they do not have a $g^2_3$ over $\bQ$.
\end{exmp}

\section{Rational Normal Scrolls under base extension}
\label{sec:rns}

Recall the following well known construction of rational
subvarieties of projective spaces (see \cite{EH87}, \cite{Xam81} or
\cite[chapter 2]{Rei97}). Let $0\le a_1\le \cdots\le a_d$ be a list
of $d$ integers with $a_d>0$ for $d\ge 1$. A rational normal scroll
of type $(a_1,\dots,a_d)$ defined over a field $k$ is a
$d$-dimensional subvariety $S_{(a_1,\cdots,a_d)}$ of $\bP^n$, for
$n=\sum_{i=1}^d a_i+d-1$, defined as follows: choose $d$
complementary linear subspaces $L_i \subset \bP^n$ for $i=1,\dots,d$
with $\dim(L_i)=a_i$. If $a_i\ne 0$, choose a rational normal curve
$C_i\subset L_i$ and an isomorphism $\phi_i:\bP^1\to C_i$ (if
$a_i=0$, set $C_i=L_i$ and $\phi_i$ to be the constant map). Then
\begin{equation}\label{scrolldef}
S_{(a_1,\cdots,a_d)}:=\bigcup_{t\in \bP^1}
\overline{\phi_1(t),\dots,\phi_d(t)},
\end{equation}
where $\overline{\phi_1(t),\dots,\phi_d(t)}$ denotes the linear
span of $\{\phi_1(t),\dots,\phi_d(t)\}\subset\bP^n$.

More abstractly,
$S_{(a_1,\cdots,a_d)}$ is the image of
$\bP(\cO_{\bP^1}(a_1)\oplus\dots\oplus\cO_{\bP^1}(a_d))$ in
projective $n-$space by the map (determined up to projective
equivalence) corresponding to the tautological line bundle $\cO(1)$.
There is a natural morphism
$$\bP(\cO_{\bP^1}(a_1)\oplus\dots\oplus\cO_{\bP^1}(a_d))\to \bP^1$$
which determines a rational map $\pi:S_{(a_1,\cdots,a_d)}\to \bP^1$
which we call structural map. With the description of \eqref{scrolldef},
each smooth point $p\in S_{(a_1,\cdots,a_d)}$ is mapped to the unique
$t\in \bP^1$ such that $p\in \overline{\phi_1(t),\dots,\phi_d(t)}$.
Any two rational normal scrolls of
the same type are projectively equivalent. It is also well known
that the degree of a rational normal scroll $S$ is equal to
$e:=\sum_{i=1}^d a_i=n-d+1=n-\dim S+1$, which is the smallest degree
of an irreducible non-degenerate $d$-fold in $\bP^n$. In fact, this
condition almost determines these subvarieties of $\bP^n$ over an
algebraically closed field, the other options being the cones over
the Veronese surface in $\bP^5$ and irreducible quadrics (noting
that quadrics of rank 3 or 4 are rational normal scrolls as well).

Recall also that a rational normal scroll is non-singular if and
only if $a_1>0$ (in which case the structural map is a morphism)
or $(a_1,\cdots,a_d)=(0,\dots,0,1)$ (this last case
since $S_{(0,\dots,0,1)}\cong \bP^n$).
Singular scrolls are cones over nonsingular scrolls.
A rational normal scroll of dimension 1 is
a rational normal curve; that is,
$S(a)\subset\bP^a$ is a rational normal curve of degree $a$.

Quadric scrolls are also well known varieties, easy to
describe.
Quadrics of rank 3, corresponding to the case
$(a_1,\cdots,a_d)=(0,\dots,0,2)$,
are cones over a conic if $d>1$;
their structural map is the projection from the
vertex of the cone $S(0,\dots,0,2)\to S(2)$. Quadrics of rank 4,
corresponding to the case $(a_1,\cdots,a_d)=(0,\dots,0,1,1)$,
are cones over the quadric surface $S(1,1)\cong \bP^1\times\bP^1$,
which supports two structures as scrolls, corresponding to two
structural maps which are the two projections to $\bP^1$. These are all
the cases with degree $e=2$.

Rational normal scrolls can also be characterized as the only
linearly normal varieties which contain a pencil of linear spaces of
codimension 1 (namely, the fibers of $\pi$).
Over an algebraically closed field, they can further
be characterized as the irreducible varieties determined by the
ideal of $2\times 2$ minors of a $2\times q$ matrix of linear forms.

A proof of the following result can be found in \cite{Rei97},
Chapter 2 and Appendix A.

\begin{lem} \label{canonic-scroll}
Let $S$ be a rational normal scroll. Then
$\Pic(S)=\bZ[H]\oplus \bZ[F]$, where $[F]$ is the class of a fiber of
a structural map, and $[H]$ is the class of a hyperplane section.
Moreover, the canonical class is $[K]=-d[H]+(e-2)[F]$.
\end{lem}

As a consequence, the structural map is unique whenever $e>2$.

We say that a subvariety $S\subset \bP^n$ is a potential rational
normal scroll if there exists a finite 
algebraic extension $L$ of
$k$ such that the base change of $S$ over $L$ is a rational normal
scroll over $L$, and we will say that $S$ splits in $L$. Or,
equivalently, that the base change of $S$ to an 
algebraic closure becomes a rational normal scroll. As it is the
case for rational normal scrolls, every potential rational normal
scroll is a cone over a nonsingular potential rational normal
scroll. A nonsingular potential rational normal scroll is thus a
$k$-form $S$ of the abstract variety
$\bP(\cO_{\bP^1}(a_1)\oplus\dots\oplus\cO_{\bP^1}(a_d))$, together
with a line bundle $L$ on $S$ which is a $k$-form of the
tautological line bundle $\cO(1)$.

\begin{exmp}
The simplest examples of potential rational normal scrolls arise in degree 2:
every quadric of rank 3 or 4 is a potential rational normal scroll.
Let us show examples of quadric potential rational normal scrolls
defined over some field $k$ which do not split over $k$.
Consider $k$ a field such that
there exists some genus 0 curve $C$ without rational points.   Then $C$ is isomorphic to a conic
in $\bP^2$, which we will denote also by $C$, which is a potential
rational normal scroll of dimension 1 (a quadric of rank 3).

From this we construct an example of dimension 2 and rank 4. Put $\bP^2$ inside
$\bP^5$ in two distinct and complementary ways, giving $L_1$ and
$L_2$ linear subspaces. If the conic is given by an equation
$ax_0^2+bx_1^2+cx_2^2=0$, then the surface is given by
$$
S: \left. \begin{array}{rl} ax_0^2+bx_1^2+cx_2^2&=0 \\
ax_3^2+bx_4^2+cx_5^2 & =0\\
x_1x_3 - x_0x_4 &=0 \\
x_1x_5 - x_2x_4 &=0\\
x_0x_5 - x_2x_3 &=0\\
ax_0x_3+bx_1x_4+cx_2x_5 & =0
\end{array}\right\} \subset \bP^5$$
\end{exmp}

\begin{thm}\label{oddscroll} Let  $S\subset \bP^n$ be a potential
rational normal scroll 
of degree $e$.
\begin{enumerate}
 \item If $S$ is a quadric of rank 4, and $\cha k \ne 2$,
there is a quadratic or a biquadratic extension $L/k$ such that
$S_L$ is a rational normal scroll, and the following are equivalent:
\begin{enumerate}
\item $S$ is a rational normal scroll.
\item $S$ is a cone over a ruled quadric surface in $\bP^3$.
\item $S$ contains a linear subspace of codimension 1.
\end{enumerate}
\item If $S$ is not a quadric of rank 4,
there is a Galois degree 2 extension $L/k$ such that $S_L$
is a rational normal scroll, and the following are
equivalent:
\begin{enumerate}
\item \label{rns} $S$ is a rational normal scroll.
\item \label{sbsp} $S$ contains a linear subspace of codimension 1.
\item \label{pt} $S$ has a $k$-rational nonsingular point.
\end{enumerate}
\end{enumerate}
Moreover,  if  $e$ is odd then $S$ is a rational normal scroll.
\end{thm}

\begin{exmp}\label{nonseparablescroll}
The hypothesis on the characteristic in the first part of the theorem can
not be dropped: let $k=\bF_2(s)$ be the
field of rational functions over the finite field $\bF_2$. The
quadric $xy+z^2+st^2=0$ in $\bP^3$ is a potential scroll which is
not a scroll over $\ksep$.
\end{exmp}

\begin{pf}
The quadratic cases are classical, but we will give some
indications. If $e=2$, then the rank is $3$ or $4$. The rank $3$
case corresponds to cones over conics, and they are scrolls if and
only if the conic has a point. The assertions in (2) are then easy.
If the rank is 4, then it is a cone over a quadric in $\bP^3$. If
the characteristic is not 2, the quadric can be diagonalized, so we
can suppose it is given by an equation of the form
$$a_0X_0^2+a_1X_1^2+a_2X_2^2+a_3X_3^2$$ for some $a_i\in K$ for $0\le
i\le 3$ with $a_0a_1a_2a_3\ne 0$. The quadric is an scroll if and
only if a change of variables can be made to get the equation
$Y_0Y_1-Y_2Y_3$, and this can be done for example over the quadratic
or biquadratic extension $L=K(\sqrt{a_0a_1},\sqrt{a_2a_3})$.

Assume that $e>2$ and $S$ is nonsingular (which is not restrictive,
as $S$ is always a cone over a nonsingular potential rational normal
scroll).

Consider the divisor class $[R']=[K]+d[H]$.
By lemma \ref{canonic-scroll},
after base change to a field $L$ over which $S$ is a rational normal
scroll, $[R']$ equals $(e-2)$ times the class of a fiber of
the structural map $\pi$. Therefore $\dim H^0(S,\cO_S(R'))=e-1$,
and the linear system $|R'|$
determines a morphism $S\rightarrow \bP^{e-2}$
whose image is a genus zero curve $D$
(which becomes a rational normal curve over $L$).

If $L/k$ is any extension where
$D$ has points, the fiber over any $L$-point of $D$ is
a linear space of codimension 1 in $S_L$; thus for every degree 2
extension where $D$ has points, $S_L$ is
a rational normal scroll, and there are
such extensions $L/k$ which are Galois.

If $S(k)\ne \emptyset$ then obviously $D$ has $k$-rational points
in $f(S(k))$, so $S$ is a rational normal scroll, which proves
$\eqref{pt}\Rightarrow\eqref{rns}$; and obviously,
$\eqref{rns}\Rightarrow\eqref{sbsp}\Rightarrow\eqref{pt}$.

Finally we prove that if $e$ is odd then $S$ has a $k$-rational point,
by induction on the dimension $d$ of $S$.
The case of dimension $d=1$ is well known, but we give here a short
argument for completeness. In this case $S$ is a curve in $\bP^e$,
which is projectively equivalent to a rational normal curve in the
algebraic closure. But then any hyperplane section of $S$ is a
$k$-rational divisor of odd degree $e$. Since $S$ is a genus $0$
curve with a divisor of odd degree, Riemann-Roch tells us that that
$S$ is isomorphic to $\bP^1$.

Now, in order to do the induction, observe first that, if $k$ is finite,
then $D$ has points, so it is not restrictive to assume $k$ infinite.
If $H$ is a hyperplane such that $H_L$ does not contain any ruling of the scroll
$S_L$, then $H_L\cap S_L$ is a scroll and $H\cap S$ a potential scroll
of the same degree $e$ and dimension $d-1$, which by the induction
hypothesis contains a $k$-point. So we are reduced to showing that
such a hyperplane exists. Now, these rulings have dimension $d-1$ so for
each of them there is a $(n-d)$-dimensional linear family of hyperplanes
containing it, and there is a $(n-d+1)$-dimensional closed subset of
the dual space $(\bP^n_L)^*$ consisting of the hyperplanes containing
\emph{some} ruling. Since $d\ge 2$, $n-d+1<n$ and by the infiniteness of $k$,
not all hyperplanes
in $(\bP^n_k)^*$ belong to this closed subset, so we are done.
\end{pf}

\section{Gonality and Rational Normal Scrolls}

Given a linearly normal projective curve $C \subset
\mathbb{P}^{g-1}$ and a map $f:C\rightarrow \bP^1$ of degree $d$,
there is a classical construction of a rational normal scroll $S
\subset \mathbb{P}^{g-1}$ containing $C$, such that $f$ is induced
by the structural pencil of $S$ (see \cite{Sch} for a detailed
exposition). The codimension of $S$ is
$\delta=h^0(C,\mathcal{O}_C(1)\otimes f^*(\mathcal
O_{\bP^1}(-1)))-1$, and it can be described as
 $$S:=\bigcup_{\lambda \in
\bP^1} \overline{f^{-1}(\lambda)} \subset \bP^{g-1},$$
where $\overline{f^{-1}(\lambda)}$ denotes the
linear span of the divisor $f^{-1}(\lambda)\subset C$ seen
as a subscheme of $\bP^{g-1}$.
So, if $C$ is canonically embedded in $\bP^{g-1}$ and
$d=\ogamma$ is the geometric gonality, then $S$ has dimension
$\ogamma-1$.
We want to show that this generalizes
for maps to genus zero curves, as follows:

\begin{prop}\label{potscroll}
Given a linearly normal projective curve $C \subset
\mathbb{P}^{g-1}$ and a map $f: C \rightarrow D$ with $D$ a genus
zero curve, there is a potential rational normal scroll $S\subset
\mathbb{P}^{g-1}$ defined over $k$, with a map $\bar f:S\rightarrow
D$ extending $f$, such that, for every extension $L/k$ with
$D_L\cong \bP^1$,
\begin{enumerate}
\item $S_L$ is a scroll, whose structural pencil is $f_L:S_L\rightarrow D_L$.
 \item The codimension of $S$ is
 $\delta=h^0(C_L,\mathcal{O}_{C_L}(1)\otimes f_L^*(\mathcal O_{\bP_L^1}(-1)))-1$.
\end{enumerate}
\end{prop}

\begin{pf}
Let $L/k$ be a Galois degree 2 extension such that $D_L\cong
\bP^1_L$, and let $S_L\subset \bP^n_L$ be the corresponding scroll
containing $C_L$. We claim that the homogeneous ideal of $S_L$ is
invariant by the action of $\Gal(L/k)$; therefore it can be
generated over $k$, defining a potential rational normal scroll
which satisfies the conditions.

To prove the claim, let us recall the construction of the ideal of $S_L$,
following \cite{EH87}. Consider the divisors of $A=f_L^{*}(p)$
where $p\in \bP^1_L$ is a point, and $B=C\cdot H-A$, where
$C\cdot H$ is a hyperplane section, on $C_L$.
Denote $V\subset H^0(\cO_{C_L}(A))$  the 2-dimensional subspace
whose projectivization is the pencil of fibers of $f$,
and $W=H^0(\cO_{C_L}(B))$. Note that $\dim W=\delta+1$.
There is a natural morphism
$\phi:V\otimes W\rightarrow H^0(\cO_{\bP^n_L}(1))$
which, after a choice of bases, can be described by
a $2 \times (\delta+1)$ matrix $M$ of linear forms on $\bP^n$.

The ideal of $S_L$ is generated by the $2\times 2$ minors of $M$,
and it is of course independent on the choice of bases.
We need to show that for every such minor, its conjugate
also belongs to the ideal.

Conjugation leaves the pencil of fibers $\bP(V)$ invariant, and so
acts on it. It follows that it also leaves the linear series
$|B|=\bP(W)$ invariant, and so acts on it. To be precise,
conjugation gives an isomorphism
$\sigma:H^0(\cO_{C_L}(B))\rightarrow H^0(\cO_{C_L}(\bar B))$, and by
invariance there is an isomorphism $i:H^0(\cO_{C_L}(\bar
B))\rightarrow H^0(\cO_{C_L}(B)))$
(product with a rational function
$f$ with $\operatorname{div}(f)=\bar B-B$). Then $\sigma_W=i\circ\sigma:
W\rightarrow W$ is an isomorphism which induces the conjugation
action on $|B|=\bP(W)$. ($\sigma_W^2$ is not in general the identity
on $W$, but multiplication by a scalar, so it does induce the
identity on $|B|$). Similarly, there is an isomorphism
$\sigma_V:V\rightarrow V$ inducing conjugation on the pencil, and if
the rational functions giving $\bar B\sim B$ and $\bar A\sim A$ are
chosen adequately, then for every $v\in V, w\in W$, with
$\phi(\sigma_V(v)\otimes \sigma_V(w))= \sigma(\phi(v\otimes w))$.
(Other choices just give proportional images).

Given bases on $V$ and $W$, application of $\sigma_V$
and $\sigma_W$ produces new bases and a new matrix
$M'$ whose minors equal the minors of $M$ conjugated.
Thus, these conjugated minors belong
to the ideal of $S_L$ as claimed.
\end{pf}

\medskip

Now we can give the second proof for theorem \ref{mainthm}
(\ref{mainthm3}).

\begin{pfm}
We know by theorem \ref{mainthm} that the gonal map factorizes
through a $k$-defined map $C\rightarrow D$ to a genus zero curve.
Consider the canonical embedding $C\subset\bP^{g-1}$, and the
potential rational normal scroll $S\subset\bP^{g-1}$ given by
proposition \ref{potscroll}. It has dimension $\ogamma-1$ and degree
$(g-1)-(\ogamma-1)+1$, so it is of odd degree, and by theorem
\ref{oddscroll}, $S$ is actually a scroll. Therefore the gonal map,
which is just the restriction to $C$ of the structural map
$S\rightarrow D$, is actually defined over $k$.
\end{pfm}

\smallskip

The geometry of the scroll $S$ allows to detect other cases in which
the gonality of $C$ must agree with its geometric gonality. We
illustrate this in the case of geometrically trigonal curves. In
this particular situation, the potential rational normal scroll is
$S\cong_{\kbar} S_{(a_1,a_2)}$, either isomorphic over $\kbar$ to
the Hirzebruch surface $\mathbb{F}_a$ if $a_1>0$, where
$a=a_2-a_1\ge0$ is an integer called the \emph{Maroni invariant} of
$C$, or isomorphic over $k$ to a (singular) cone over a
genus zero curve if
$a_1=0$. Note that in the geometrically trigonal case the potential
scroll is defined by the degree 2 part of the ideal of $C$ in
$\bP^{n-1}$, so it is defined over $k$ even in the case proposition
\ref{potscroll} cannot not be applied.

\begin{thm}
\label{theoremtrig}
 Let $C$ be a geometrically trigonal (i.e., such that $\ggamma=3$)
 curve with genus $g>4$. Then its gonality $\gamma$ satisfies
 $\gamma\le 6$ and if either
\begin{enumerate}
\item $g$ is odd, or
\item $g$ is even, and the Maroni invariant $a$ of $C$ satisfies
$a>0$ and $g+a\in (4)$,
\end{enumerate}
then $\gamma=3$.
\end{thm}

\begin{pf}  Note that the degree of $S$ equals $e=a_1+a_2=g-2>2$, so
$S$ is already a scroll over $\ksep$ by theorem \ref{oddscroll}.

We consider first the singular case $a_1=0$. In this case
$\ogamma=3$, hence $\gamma\le 6$ by corollary \ref{Conicsplitting}
and if $g$ is odd, $\gamma=3$ by proposition \ref{mainthm2} . The
case (2) cannot occur since $g+a_2=2g-2\equiv 2 (4)$.

 Suppose now $S$ is non-singular. Therefore $S$ is isomorphic
 over $\ksep$ to $\mathbb{F}_a$.
 Let $F$ be a fiber of the ruling on $S_{\ksep}$, defined over
 $\ksep$.
%
We want to show that under the
hypotheses, $F$ is defined over $k$ and so the gonal map (which is
the restriction of the map given by $|F|$ to $C$) is defined  over
$k$. In fact it is enough to see that an odd multiple $rF$ is
defined over $k$; indeed,  in that case one obtains a map
$S\rightarrow \mathbb{P}^r$ whose image is the rational normal curve
which, if $r$ is odd, is isomorphic to  $\mathbb{P}^1$.

Denote by $H$ a hyperplane section, by $K$ a canonical divisor on
$S\cong \mathbb{F}_a$, and by $E$ the negative section
($E^2=-a<0$, which is the only smooth curve in $\mathbb{F}_a$ with
negative selfintersection and is therefore Galois-invariant);
all three are divisors defined over $k$.
By Lemma \ref{canonic-scroll} we
have over $\ksep$ that $\operatorname{Pic}(S_{\ksep})\cong
\mathbb{Z}H\oplus\mathbb{Z}F$. In fact, it is well known
(see e.g. \cite[V.2]{Har} or \cite{Rei97})
that $H$ is
linearly equivalent to $E+a_2 F$, whereas $K$ is linearly equivalent
to $-2E-(a+2)F$. Thus, whenever $a_2$ or $a$ are
odd (which is satisfied under the given hypotheses) there is an odd
multiple of $F$ defined over $k$.
\end{pf}

\smallskip
If $g=4$, the bound for the gonality $\gamma$ also holds since
$2g-2=6$ in this case; but if $S$ is non-singular, the Maroni
invariant is $a=0$ although the gonal map is not unique and in this case
the gonality over $k$ can be $6$; in fact the scroll $S$ belongs to the
case considered in theorem \ref{oddscroll}, (1).

\end{document}